\documentclass[12pt,leqno]{article}
\usepackage{amsmath, amsthm, amsfonts, amssymb, color}
\usepackage{mathrsfs}
\theoremstyle{definition}

\newcommand{\scr}[1]{\mathscr #1}
\definecolor{wco}{rgb}{0.5,0.2,0.3}

\numberwithin{equation}{section} \theoremstyle{remark}

\newcommand{\ua}{\uparrow}

\title{{\bf A Note on the Paper $``$Poincar\'e Inequality on the Path Space of Poisson Point
Processes''}}

\author{
{\bf Feng-Yu Wang  and Chenggui Yuan}\\
\footnotesize{   Department of Mathematics,
Swansea University, Singleton Park, SA2 8PP, UK}\\
\footnotesize{Email:
F.Y.Wang@swansea.ac.uk; \ C.Yuan@swansea.ac.uk}
}\date{\ }
\begin{document}
\def\R{\mathbb R}  \def\ff{\frac} \def\ss{\sqrt} \def\B{\mathbf
B}
\def\N{\mathbb N} \def\kk{\kappa} \def\m{{\bf m}}
\def\dd{\delta} \def\DD{\Delta} \def\vv{\varepsilon} \def\rr{\rho}
\def\<{\langle} \def\>{\rangle} \def\GG{\Gamma} \def\gg{\gamma}
  \def\nn{\nabla} \def\pp{\partial} \def\EE{\scr E}
\def\d{\text{\rm{d}}} \def\bb{\beta} \def\aa{\alpha} \def\D{\scr D}
  \def\si{\sigma} \def\ess{\text{\rm{ess}}}
\def\beg{\begin} \def\beq{\begin{equation}}  \def\F{\scr F}
\def\Ric{\text{\rm{Ric}}} \def\Hess{\text{\rm{Hess}}}
\def\e{\text{\rm{e}}} \def\ua{\underline a} \def\OO{\Omega}  \def\oo{\omega}
 \def\tt{\tilde} \def\Ric{\text{\rm{Ric}}}
\def\cut{\text{\rm{cut}}} \def\P{\mathbb P} \def\ifn{I_n(f^{\bigotimes n})}
\def\C{\mathbb C}      \def\aaa{\mathbf{r}}     \def\r{r}
\def\gap{\text{\rm{gap}}} \def\prr{\pi_{{\bf m},\varrho}}  \def\r{\mathbf r}
\def\Z{\mathbb Z} \def\vrr{\varrho} \def\ll{\lambda}
\def\L{\scr L}\def\Tt{\tt} \def\TT{\tt}\def\II{\mathbb I}
\def\i{{\rm in}}\def\Sect{{\rm Sect}}\def\E{\mathbb E} \def\H{\mathbb H}
\def\M{\scr M}\def\Q{\mathbb Q} \def\texto{\text{o}} \def\LL{\Lambda}
\def\Rank{{\rm Rank}} \def\B{\scr B}
\def\T{\mathbb T}\def\i{{\rm i}} \def\ZZ{\hat\Z}

\maketitle


 In the recent paper \cite{WY} we proved a Poincar\'e inequality on the path space of   compound
 Poisson processes by using transition probabilities and the Markov property. Our purpose is to develop a general argument for stochastic analysis on the path space of jump processes as indicated in the introduction.
 
 After the paper was published, we learnt from some colleagues that the Poincar\'e inequality
  we obtained was already known by L. Wu \cite[Remark 1.4]{W}
   using martingale representations of Poisson point processes.
   Moreover, the invalidity of the log-Sobolev inequality was already known by D. Surgailis \cite{S}.
    The Dirichlet form considered in these two papers is the birth-death Dirichlet form
    on the $L^2$ space of a Poisson measure, which covers our framework by taking
    the Poisson measure with intensity   $\nu(\d x)\times \d t$ on $\R^d\times [0,T]$,
    where $\nu$ is the L\'evy measure of the underlying compound Poisson process.
   Unfortunately, we did not aware this point when we prepared our paper, so that the important references \cite{W,S} were not cited.

 \paragraph{Acknowledgement} We would  like to thank Professors E. Lytvynov, M. R\"ockner and G. Last for their kindly comments.

 \footnotesize
\beg{thebibliography}{99}

\bibitem{S}
D. Surgailis, \emph{On the multiple Poisson stochastic integrals and associated Markov semigroups},
Proba. and Math. Sta. 3 (1984), 217-239.

\bibitem{WY}  F.-Y. Wang, C. Yuan, \emph{Poincar\'e inequality on the path space of Poisson point
processes,} 23(2010), 824--833.

\bibitem{W}
L. Wu, \emph{A new modified logarithmic Sobolev inequality for
Poisson point processes and several applications}, Probab. Theory
Relat. Fields 118 (2000), 427-438.

 \end{thebibliography}
\end{document}